\input amstex
\input amsppt.sty
\magnification=\magstep1
\hsize=32truecc
\vsize=22.2truecm
\baselineskip=16truept
\NoBlackBoxes
\TagsOnRight \pageno=1 \nologo
\def\Z{\Bbb Z}

\def\l{\left}
\def\r{\right}
\def\bg{\bigg}
\def\({\bg(}
\def\[{\bg\lfloor}
\def\){\bg)}
\def\]{\bg\rfloor}
\def\t{\text}
\def\f{\frac}

\def\bi{\binom}
\def\eq{\equiv}

\def\ls{\leqslant}

\def\mo{\roman{mod}}

\def\Proof{\noindent{\it Proof}}

\def\Remark{\medskip\noindent{\it  Remark}}

\def\Ack{\medskip\noindent {\bf Acknowledgment}}
\hbox {Adv. in Appl. Math. 51(2013), no.\,4, 524--535.}
\bigskip
\topmatter
\title Congruences for Franel numbers \endtitle
\author Zhi-Wei Sun\endauthor
\leftheadtext{Zhi-Wei Sun} \rightheadtext{Congruences for Franel
numbers}
\affil Department of Mathematics, Nanjing University\\
 Nanjing 210093, People's Republic of China
  \\  zwsun\@nju.edu.cn
  \\ {\tt http://math.nju.edu.cn/$\sim$zwsun}
\endaffil
\abstract The Franel numbers given by $f_n=\sum_{k=0}^n\bi nk^3$
 $(n=0,1,2,\ldots)$ play important roles in both combinatorics and number theory.
  In this paper we initiate the systematic investigation of fundamental congruences
 for the Franel numbers. We mainly establish for any prime $p>3$ the following
 congruences:
 $$\gather \sum_{k=0}^{p-1}(-1)^kf_k\eq\l(\f p3\r)\pmod{p^2},
\ \ \sum_{k=0}^{p-1}(-1)^kkf_k\eq-\f23\l(\f p3\r)\pmod{p^2},
\\\sum_{k=1}^{p-1}\f{(-1)^k}kf_k\eq0\pmod {p^2},
\ \ \sum_{k=1}^{p-1}\f{(-1)^k}{k^2}f_k\eq0\ (\mo\ p).
\endgather$$
\endabstract
\thanks 2010 {\it Mathematics Subject Classification}. \,Primary 11A07, 11B65;
Secondary  05A10, 11B37, 11B75.
\newline\indent Supported by the National Natural Science
Foundation (grant 11171140) of China and the PAPD of Jiangsu Higher Education Institutions.
\endthanks

\endtopmatter
\document

\heading{1. Introduction}\endheading

In 1894, Franel [F] noted that the numbers
$$f_n=\sum_{k=0}^n\bi nk^3\ \ \ (n=0,1,2,\ldots)\tag1.1$$
(cf. [Sl, A000172]) satisfy the recurrence relation:
$$(n+1)^2f_{n+1}=(7n^2+7n+2)f_n+8n^2f_{n-1}\ \ (n=1,2,3,\ldots).\tag1.2$$
Such numbers are now called Franel numbers. For a combinatorial
interpretation of the Franel numbers, see Callan [C]. Recall
that the Ap\'ery numbers given by
$$A_n=\sum_{k=0}^n\bi nk^2\bi{n+k}k^2=\sum_{k=0}^n\bi{n+k}{2k}^2\bi{2k}k^2\
(n=0,1,2,\ldots)$$ were introduced by Ap\'ery [A], and
they can be expressed in terms of Franel numbers as follows:
$$A_n=\sum_{k=0}^n\bi nk\bi{n+k}kf_k\tag1.3$$
(see Strehl [St92]). The Franel numbers are also related to the theory of modular forms,
see, e.g., Zagier [Z].

In this paper we study congruences for the Franel
numbers systematically. As usual, for any odd prime $p$ and integer $a$,
$(\f ap)$ denotes the Legendre symbol, and $q_p(a)$ stands for the Fermat quotient $(a^{p-1}-1)/p$
if $p\nmid a$.

Now we state our main result.

\proclaim{Theorem 1.1} Let $p>3$ be a prime. For any $p$-adic integer $r$ we have
$$\sum_{k=0}^{p-1}(-1)^k\bi{k+r}kf_k\eq\sum_{k=0}^{p-1}\bi{2k}k\bi{k+r}k^2\pmod{p^2}.\tag1.4$$
In particular,
$$\align \sum_{k=0}^{p-1}(-1)^kf_k\eq&\l(\f p3\r)\pmod{p^2},\tag1.5
\\\sum_{k=0}^{p-1}(-1)^kkf_k\eq&-\f23\l(\f p3\r)\pmod{p^2},\tag1.6
\\\sum_{k=0}^{p-1}(-1)^kk^2f_k\eq&\f{10}{27}\l(\f p3\r)\pmod{p^2},\tag1.7
\endalign$$
and
$$\sum_{k=0}^{p-1}\f{\bi{2k}kf_k}{(-4)^k}\eq\sum_{k=0}^{p-1}\f{\bi{2k}k^3}{16^k}\pmod{p^2}.\tag1.8$$
We also have
 $$\align\sum_{k=1}^{p-1}\f{(-1)^k}kf_k\eq&0\pmod {p^2},\tag1.9
\\\sum_{k=1}^{p-1}\f{(-1)^k}{k^2}f_k\eq&0\pmod {p},\tag1.10
\\\sum_{k=1}^{p-1}\f{(-1)^k}kf_{k-1}\eq&3q_p(2)+3p\,q_p(2)^2\pmod{p^2},\tag1.11
\endalign$$
\endproclaim
\Remark\ 1.1. Fix a prime $p>3$. In contrast with (1.5), we conjecture that
$$\sum_{n=0}^{p-1}(-1)^n\sum_{k=0}^n\bi nk^3(-8)^k\eq \sum_{k=0}^{p-1}\f{f_k}{8^k} \eq\l(\f p3\r)\pmod{p^2}.$$
As $f_k\eq(-8)^kf_{p-1-k}\ (\mo\ p)$ for all $k=0,\ldots,p-1$ by [JV, Lemma 2.6],
(1.11) implies that
$$\sum_{k=1}^{p-1}\f{f_k}{k8^k}\eq\sum_{k=1}^{p-1}\f{(-1)^k}kf_{p-1-k}=\sum_{k=1}^{p-1}\f{(-1)^{p-k}}{p-k}f_{k-1}\eq 3q_p(2)\pmod p.$$
Motivated by (1.5) and (1.6), we conjecture that both $(\sum_{k=0}^{n-1}(-1)^kf_k)/n^2$ and
$(\sum_{k=0}^{n-1}(-1)^kkf_k)/n^2$ are $3$-adic integers for any positive integer $n$.
Concerning (1.8) the author [S11, Conj. 5.2(ii)] conjectured that
$$\sum_{k=0}^{p-1}\f{\bi{2k}k^3}{16^k}\eq\cases 4x^2-2p\ (\mo\ p^2)&\t{if}\ p=x^2+3y^2\ (x,y\in\Z),
\\0\ (\mo\ p^2)&\t{if}\ p\eq 2\ (\mo\ 3).\endcases$$
See also [S13] for other connections between $p=x^2+3y^2$ and Franel numbers.
(1.10) can be extended as
$$\sum_{k=1}^{p-1}\f{(-1)^{kr}}{k^{r-1}}f_k^{(r)}\eq0\pmod p,\tag1.12$$
where $r$ is any positive integer and $f_k^{(r)}:=\sum_{j=0}^k\bi kj^r$. Note that $f_k^{(2)}=\bi{2k}k$
and $\sum_{k=1}^{p-1}\bi{2k}k/k\eq0\ (\mo\ p^2)$ by [ST10].
\medskip

Let $p>3$ be a prime.
Similar to (1.5)-(1.7), we are also able to show that
$$\sum_{k=0}^{p-1}(-1)^kk^3f_k\eq-\f{10}{81}\l(\f p3\r)\ (\mo\ p^2)
\ \t{and}\ \sum_{k=0}^{p-1}(-1)^kk^4f_k\eq-\f{14}{243}\l(\f p3\r)\ (\mo\ p^2).$$
In general, for any positive integer $r$ and prime $p>\max\{r,3\}$ there should be an odd integer $a_r$ (not dependent on $p$) such that
$$\sum_{k=0}^{p-1}(-1)^kk^rf_k\eq\f{2a_r}{3^{2r-1}}\l(\f p3\r)\pmod{p^2}.$$

\heading{2. Proof of Theorem 1.1}\endheading

We first establish an auxiliary theorem on the polynomials
$$f_n(x):=\sum_{k=0}^n\bi nk^2\bi{2k}nx^k=\sum_{k=0}^n\bi nk\bi k{n-k}\bi{2k}kx^k\ \ (n=0,1,2,\ldots).$$

\proclaim{Theorem 2.1} Let $p$ be an odd prime and let $r$ be any $p$-adic integer. Then
$$\sum_{l=0}^{p-1}(-1)^l\bi{l+r}lf_l(x)\eq\sum_{k=0}^{p-1}\bi{2k}kx^k\bi{k+r}k^2\pmod{p^2}.\tag2.1$$
\endproclaim
\Proof. Observe that
$$\align\sum_{l=0}^{p-1}(-1)^l\bi{l+r}lf_l(x)=&\sum_{l=0}^{p-1}(-1)^l\bi{l+r}l\sum_{k=0}^l\bi lk\bi k{l-k}\bi{2k}kx^k
\\=&\sum_{k=0}^{p-1}\bi{2k}kx^k\sum_{l=k}^{\min\{2k,p-1\}}(-1)^l\bi lk\bi{k}{l-k}\bi{l+r}l.
\endalign$$
If $(p-1)/2<k\ls p-1$ and $p\ls l\ls 2k$, then
$$\bi{2k}k=\f{(2k)!}{(k!)^2}\eq0\ (\mo\ p)
\ \ \t{and}\ \ \bi{l}k=\f{l!}{k!(l-k)!}\eq0\ (\mo\ p).$$
Thus
$$\sum_{l=0}^{p-1}(-1)^l\bi{l+r}lf_l(x)\eq\sum_{k=0}^{p-1}\bi{2k}kx^k\sum_{l=k}^{2k}(-1)^l\bi lk\bi{k}{l-k}\bi{l+r}l\pmod{p^2},$$
and hence it suffices to show the identity
$$\sum_{l=k}^{2k}(-1)^l\bi lk\bi k{l-k}\bi{x+l}l=\bi{x+k}k^2.\tag2.2$$

By the well-known Chu-Vandermonde identity (cf. (3.1) of [G, p.22]),
$$\sum_{j=0}^k\bi{y}j\bi{z}{k-j}=\bi{y+z}k.$$
Therefore
$$\align &\sum_{l=k}^{2k}(-1)^l\bi lk\bi k{l-k}\bi{x+l}l
\\=&\sum_{l=k}^{2k}\bi lk\bi{k}{l-k}\bi{-x-1}l
=\bi{-x-1}k\sum_{l=k}^{2k}\bi{-x-1-k}{l-k}\bi k{l-k}
\\=&\bi{-x-1}k\sum_{j=0}^k\bi{-x-1-k}{j}\bi k{k-j}=\bi{-x-1}k^2=\bi{x+k}k^2.
\endalign$$
This proves (2.2) and hence (2.1) follows. \qed

\proclaim{Lemma 2.1} For any nonnegative integer $n$, the integer $f_n(1)$ coincides with the Franel number $f_n$.
\endproclaim
\Proof. The identity $\sum_{k=0}^n\bi nk^2\bi{2k}n=f_n$ is a known result due to Strehl [St94]. \qed

\proclaim{Lemma 2.2} For each positive integer $m$ we have
$$\sum_{k=0}^{n-1}P_m(k)\bi{2k}k=n^m\bi{2n}n\quad\t{for all}\ n=1,2,3,\ldots,$$
where $P_m(x):=2(2x+1)(x+1)^{m-1}-x^m$.
\endproclaim
\Proof. The desired result follows immediately by induction on $n$. \qed

\proclaim{Lemma 2.3} Let $m$ be a positive integer. For $n=0,1,\ldots,m$ we have
$$\sum_{k=0}^n\bi xk\bi{-x}{m-k}=\f{m-n}m\bi{x-1}n\bi{-x}{m-n}.$$
\endproclaim
\Remark\ 2.1. This is a known result due to Andersen, see, e.g., (3.14) of [G, p.\,23].

\proclaim{Lemma 2.4 {\rm ([S11, Lemma 2.1])}} Let $p$ be an odd
prime. For any $k=1,\ldots,p-1$ we have
$$k\bi{2k}k\bi{2(p-k)}{p-k}\eq(-1)^{\lfloor 2k/p\rfloor-1}2p\pmod{p^2}.$$
\endproclaim

Recall that the harmonic numbers and the second-order harmonic numbers are given by
$$H_n=\sum_{0<k\ls n}\f1k\ \ \t{and}\ \ H_n^{(2)}=\sum_{0<k\ls n}\f1{k^2}\ \ \ (n=0,1,2,\ldots)$$
respectively. Let $p>3$ be a prime. In 1862, Wolstenholme [W] proved that
$$H_{p-1}\eq0\ (\mo\ p^2)\ \ \t{and}\ \ H_{p-1}^{(2)}\eq0\ (\mo\ p).$$
Note that
$$H_{(p-1)/2}^{(2)}\eq\f12\sum_{k=1}^{(p-1)/2}\l(\f1{k^2}+\f1{(p-k)^2}\r)=\f12H_{p-1}^{(2)}\eq0\pmod p.$$
In 1938, Lehmer [L] showed that
$$H_{(p-1)/2}\eq -2q_p(2)+p\,q_p(2)^2\pmod{p^2}.\tag2.3$$

\proclaim{Lemma 2.5} Let $p>3$ be a prime. Then
$$f_{p-1}\eq1+3p\,q_p(2)+3p^2q_p(2)^2\pmod{p^3}.\tag2.4$$
\endproclaim
\Proof. For any $k=1,\ldots,p-1$, we obviously have
$$\align&(-1)^k\bi{p-1}k=\prod_{j=1}^k\l(1-\f pj\r)
\\\eq& 1-pH_k+\f {p^2}2\sum_{1\ls i<j\ls k}\f2{ij}
=1-pH_k+\f{p^2}2\l(H_k^2-H_k^{(2)}\r)\pmod{p^3}.
\endalign$$
Thus
$$\align f_{p-1}-1=&\sum_{k=1}^{p-1}\bi{p-1}k^3\eq\sum_{k=1}^{p-1}(-1)^k\l(1-pH_k+\f{p^2}2\l(H_k^2-H_{k}^{(2)}\r)\r)^3
\\\eq&-3p\sum_{k=1}^{p-1}(-1)^kH_k+\f{9}2p^2\sum_{k=1}^{p-1}(-1)^kH_k^2-\f 32p^2\sum_{k=1}^{p-1}(-1)^kH_k^{(2)}\pmod{p^3}.
\endalign$$
Clearly
$$\align \sum_{k=1}^{p-1}(-1)^kH_k=&\sum_{k=1}^{p-1}\sum_{j=1}^k\f{(-1)^k}j
=\sum_{j=1}^{p-1}\f{\sum_{k=j}^{p-1}(-1)^k}j=\sum^{p-1}\Sb j=1\\2\mid j\endSb\f1j
\\=&\f12H_{(p-1)/2}\eq-q_p(2)+\f p2q_p(2)^2\pmod{p^2}\ \ \t{(by (2.3))}
\endalign$$
and
$$\sum_{k=1}^{p-1}(-1)^kH_k^{(2)}=\sum_{j=1}^{p-1}\f{\sum_{k=j}^{p-1}(-1)^k}{j^2}=\sum_{i=1}^{(p-1)/2}\f1{(2i)^2}=\f{H_{(p-1)/2}^{(2)}}4\eq0\pmod p.$$
Observe that
$$\align\sum_{k=1}^{p-1}(-1)^kH_k^2=&\sum_{k=1}^{p-1}(-1)^{p-k}H_{p-k}^2=\sum_{k=1}^{p-1}(-1)^{k-1}\(H_{p-1}-\sum_{0<j<k}\f1{p-j}\)^2
\\\eq&-\sum_{k=1}^{p-1}(-1)^k\l(H_k-\f1k\r)^2
\\=&-\sum_{k=1}^{p-1}(-1)^kH_k^2+2\sum_{k=1}^{p-1}\f{(-1)^k}kH_k-\sum_{k=1}^{p-1}\f{(-1)^k}{k^2}\pmod{p}.
\endalign$$
Clearly,
$$\sum_{k=1}^{p-1}\f{(-1)^k}{k^2}\eq\sum_{k=1}^{p-1}\f{1+(-1)^k}{k^2}=\sum_{j=1}^{(p-1)/2}\f2{(2j)^2}\eq0\pmod p,$$
and
$$\sum_{k=1}^{p-1}\f{(-1)^k}kH_k=\sum^{p-1}\Sb k=1\\2\mid k\endSb\f{H_k}k-\sum^{p-1}\Sb k=1\\2\nmid k\endSb\f{H_k}k
\eq\f{q_p(2)^2}2-\l(-\f{q_p(2)^2}2\r)\pmod p$$
by [S12a, Lemma 2.3]. Therefore
$$\sum_{k=1}^{p-1}(-1)^kH_k^2\eq\sum_{k=1}^{p-1}\f{(-1)^k}kH_k\eq q_p(2)^2\pmod p.$$

Combining the above, we finally obtain
$$\align f_{p-1}-1\eq&-3p\l(-q_p(2)+\f p2q_p(2)^2\r)+\f{9}2p^2q_p(2)^2\pmod{p^3}
\endalign$$
and hence (2.4) holds. \qed

\proclaim{Lemma 2.6} Let $p$ be any prime. Then
$$\bi{p-1}k\bi{p+k}k\eq(-1)^k\pmod{p^2}\quad\t{for}\ k=0,1,\ldots,p-1,$$
and
$$\bi{2k}k\sum_{n=k}^{p-1}(2n+1)\bi{n+k}{2k}\eq p^2\,\f{(-1)^k}{k+1}\pmod{p^4}\quad\t{for}\ k=0,\ldots,p-2.$$
\endproclaim
\Proof. Let $k\in\{0,1,\ldots,p-1\}$. Clearly
$$\align\bi{p-1}k\bi{p+k}k=&\prod_{0<j\ls k}\l(\f{p-j}j\cdot\f{p+j}j\r)
\eq(-1)^k\pmod{p^2}.
\endalign$$
In view of the known identity
$\sum_{n=0}^m\bi nl=\bi{m+1}{l+1}$ $(l,m=0,1,\ldots)$
(see, e.g., (1.52) of [G, p.\,7]) which can be easily proved by induction, we have
$$\align\sum_{n=k}^{p-1}\f{2n+1}{2k+1}\bi{n+k}{2k}=&\sum_{n=k}^{p-1}\l(\f{2(n+k+1)}{2k+1}-1\r)\bi{n+k}{2k}
\\=&2\sum_{n=k}^{p-1}\bi{n+k+1}{2k+1}-\sum_{n=k}^{p-1}\bi{n+k}{2k}
\\=&2\bi{p+k+1}{2k+2}-\bi{p+k}{2k+1}=\f p{k+1}\bi{p+k}{2k+1}
\endalign$$
and hence
$$\bi{2k}k\sum_{n=k}^{p-1}(2n+1)\bi{n+k}{2k}=p\f{2k+1}{k+1}\bi{2k}k\bi{p+k}{2k+1}=\f {p^2}{k+1}\bi{p-1}k\bi{p+k}k.$$
Thus, if $k<p-1$ then
$$\bi{2k}k\sum_{n=k}^{p-1}(2n+1)\bi{n+k}{2k}\eq\f {p^2}{k+1}(-1)^k\pmod{p^4}$$
as desired. \qed

\medskip
\noindent{\it Proof of Theorem} 1.1. In view of Lemma 2.1,
(2.1) with $x=1$ gives (1.4).

(2.1) with $r=0$ yields the congruence
$$\sum_{k=0}^{p-1}(-1)^kf_k(x)\eq\sum_{k=0}^{p-1}\bi{2k}kx^k\pmod{p^2}.$$
In the case $x=1$, this gives (1.5) since
$\sum_{k=0}^{p-1}\bi{2k}k\eq (\f p3)\ (\mo\ p^2)$ by [ST11, (1.9)].

By (2.1) with $r=0,1$,
$$\align&\sum_{k=0}^{p-1}(3(k+1)-1)(-1)^kf_k(x)
\\\eq&\sum_{k=0}^{p-1}\bi{2k}kx^k\l(3(k+1)^2-1\r)
=\sum_{k=0}^{p-1}P_2(k)\bi{2k}kx^k\pmod{p^2}
\endalign$$
where $P_2(x)=2(2x+1)(x+1)-x^2=3x^2+6x+2$.
Thus, with the help of Lemmas 2.1-2.2, we have
$$\sum_{k=0}^{p-1}(3k+2)(-1)^kf_k\eq 0\pmod{p^2}\tag2.5$$
and hence (1.6) holds in view of (1.5).

Taking $r=2$ in (2.1) we get
$$2\sum_{k=0}^{p-1}(k^2+3k+2)(-1)^kf_k(x)\eq\sum_{k=0}^{p-1}\bi{2k}kx^k((k+1)(k+2))^2\pmod{p^2}.$$
In view of (2.5), this yields
$$2\sum_{k=0}^{p-1}(-1)^kk^2f_k\eq\sum_{k=0}^{p-1}\bi{2k}k(k^2+3k+2)^2\pmod{p^2}.$$
Note that
$$27(k^2+3k+2)^2=9P_4(k)+12P_3(k)+23P_2(k)+20$$
where $P_m(x)$ is given by Lemma 2.2. Therefore, with the help of Lemma 2.3 and [ST11, (1.9)], we have
$$54\sum_{k=0}^{p-1}(-1)^kk^2f_k\eq\sum_{k=0}^{p-1}(9P_4(k)+12P_3(k)+23P_2(k)+20)\bi{2k}k\eq 20\l(\f p3\r)\ (\mo\ p^2)$$
and hence (1.7) follows.

Putting $r=-1/2$ in (2.1) and noting that
$\bi{k-1/2}k=\bi{2k}k/4^k$,
we then obtain
$$\sum_{k=0}^{p-1}\f{\bi{2k}kf_k(x)}{(-4)^k}\eq\sum_{k=0}^{p-1}\f{\bi{2k}k^3}{16^k}x^k\pmod{p^2}.\tag2.6$$
In the case $x=1$ this gives (1.8).

Now we prove (1.9). Observe that
$$\sum_{l=1}^{p-1}\f{(-1)^l}l\sum_{k=0}^l\bi lk\bi{k}{l-k}\bi{2k}kx^k
=\sum_{k=1}^{p-1}\f{\bi{2k}k}kx^k\sum_{l=k}^{p-1}(-1)^l\bi{l-1}{k-1}\bi{k}{l-k}.$$
If $1\ls k\ls(p-1)/2$, then
$$\align\sum_{l=k}^{p-1}(-1)^l\bi{l-1}{k-1}\bi k{l-k}
=&\sum_{l=k}^{2k}(-1)^l\bi{l-1}{k-1}\bi k{l-k}
\\=&\sum_{j=0}^k(-1)^{k+j}\bi{k+j-1}j\bi kj
\\=&(-1)^k\sum_{j=0}^k\bi{-k}j\bi{k}{k-j}=(-1)^k\bi0k=0
\endalign$$
by the Chu-Vandermonde identity.
If $(p+1)/2\ls k\ls p-1$, then
$$\align\sum_{l=k}^{p-1}(-1)^l\bi{l-1}{k-1}\bi k{l-k}
=&\sum_{j=0}^{p-1-k}(-1)^{k+j}\bi{k+j-1}j\bi kj
\\=&(-1)^k\sum_{j=0}^{p-1-k}\bi{-k}j\bi{k}{k-j}
\endalign$$
and hence applying Lemma 2.3 we get
$$\align&\sum_{l=k}^{p-1}(-1)^l\bi{l-1}{k-1}\bi k{l-k}
\\=&(-1)^k\f{k-(p-1-k)}k\bi{-k-1}{p-1-k}\bi k{k-(p-1-k)}
\\=&(-1)^{p-1}\l(\f {p-k}k\r)^2\bi{p-1}{k-1}\bi k{p-k}
\\\eq&(-1)^{k-1}\bi k{p-k}=\bi{p-2k-1}{p-k}
\\\eq&\bi{2(p-k)-1}{p-k}=\f12\bi{2(p-k)}{p-k}\pmod p.
\endalign$$
Note that $\bi{2k}k\eq0\ (\mo\ p)$ for $k=(p+1)/2,\ldots,p-1$. By the above,
$$\sum_{l=1}^{p-1}\f{(-1)^l}lf_l(x)\eq\sum_{k=(p+1)/2}^{p-1}\f{\bi{2k}k}kx^k\f{\bi{2(p-k)}{p-k}}2
\eq p\sum_{k=(p+1)/2}^{p-1}\f{x^k}{k^2}\pmod{p^2}\tag2.7$$
with the help of Lemma 2.4. Hence (1.9) follows from (2.7) in the case $x=1$ since
$$2\sum_{k=(p+1)/2}^{p-1}\f1{k^2}\eq\sum_{k=(p+1)/2}^{p-1}\l(\f1{k^2}+\f1{(p-k)^2}\r)=H_{p-1}^{(2)}\eq0\pmod p.$$

Instead of proving (1.10) we show its extension (1.12). Clearly,
$$\sum_{k=1}^{p-1}\f{(-1)^{kr}}{k^{r-1}}=\sum_{k=1}^{(p-1)/2}\l(\f{(-1)^{kr}}{k^{r-1}}+\f{(-1)^{(p-k)r}}{(p-k)^{r-1}}\r)\eq0\pmod p.$$
Thus
$$\align \sum_{l=1}^{p-1}\f{(-1)^{lr}}{l^{r-1}}f_l^{(r)}\eq&\sum_{l=1}^{p-1}\f{(-1)^{lr}}{l^{r-1}}\sum_{k=1}^l\bi{l}{k}^r
=\sum_{k=1}^{p-1}\f{1}{k^{r-1}}\sum_{l=k}^{p-1}(-1)^{lr}\bi{l-1}{k-1}^{r-1}\bi lk
\\=&\sum_{k=1}^{p-1}\f{1}{k^{r-1}}\sum_{j=0}^{p-1-k}(-1)^{(k+j)r}\bi{k+j-1}j^{r-1}\bi{k+j}j
\\=&\sum_{k=1}^{p-1}\f{(-1)^{kr}}{k^{r-1}}\sum_{j=0}^{p-1-k}\bi{-k}j^{r-1}\bi{-k-1}j
\\\eq&\sum_{k=1}^{p-1}\f{(-1)^{kr}}{k^{r-1}}\sum_{j=0}^{p-k-1}\bi{p-k}j^{r-1}\bi{p-k-1}j\pmod p.
\endalign$$
For any positive integer $n$, we have
$$f_n^{(r)}=\sum_{k=0}^n\l(\f kn+\f{n-k}n\r)\bi nk^r=2\sum_{k=0}^n\f{n-k}n\bi nk^r=2\sum_{k=0}^{n-1}\bi{n}k^{r-1}\bi {n-1}k.$$
Therefore,
$$\align \sum_{l=1}^{p-1}\f{(-1)^{lr}}{l^{r-1}}f_l^{(r)}\eq&\sum_{k=1}^{p-1}\f{(-1)^{kr}}{k^{r-1}}\cdot\f{f_{p-k}^{(r)}}2
=\f12\sum_{k=1}^{p-1}\f{(-1)^{(p-k)r}f_k^{(r)}}{(p-k)^{r-1}}
\\\eq&-\f12\sum_{k=1}^{p-1}\f{(-1)^{kr}}{k^{r-1}}f_k^{(r)}\pmod p
\endalign$$
and hence (1.12) follows.

Finally we show (1.11). By (1.3) and Lemma 2.6,
$$\align\f1p\sum_{n=0}^{p-1}(2n+1)A_n=&\f1p\sum_{n=0}^{p-1}(2n+1)\sum_{k=0}^n\bi{n+k}{2k}\bi{2k}kf_k
\\=&\f1p\sum_{k=0}^{p-1}\bi{2k}kf_k\sum_{n=k}^{p-1}(2n+1)\bi{n+k}{2k}
\\\eq&\f{f_{p-1}}p\bi{2p-2}{p-1}(2p-1)+p\sum_{k=0}^{p-2}\f{(-1)^kf_k}{k+1}
\\=&\bi{2p-1}{p-1}f_{p-1}-p\sum_{k=1}^{p-1}\f {(-1)^k}kf_{k-1}\pmod{p^3}.
\endalign$$
Combining this with Wolstenholme's congruence $\bi{2p-1}{p-1}\eq1\ (\mo\ p^3)$ (cf. [W]) and [S12b, (1.6)] we obtain
$$\sum_{k=1}^{p-1}\f{(-1)^{k}f_{k-1}}{k}\eq\f{f_{p-1}-1}p\eq 3q_p(2)+3p\,q_p(2)^2\pmod{p^2}$$
by Lemma 2.5.
\qed

\Ack. The author would like to thank Prof. Qing-Hu Hou at Nankai Univ. for his comments on
the author's original version of Lemma 2.2.


\widestnumber\key{ST11}

 \Refs

\ref\key A\by R. Ap\'ery\paper {\rm Irrationalit\'e de $\zeta(2)$ et $\zeta(3)$}
\jour Ast\'erisque 61(1979) 11--13\endref

\ref \key C\by D. Callan\paper {\rm A combinatorial interpretation
for an identity of Barrucand} \jour J. Integer Seq. 11(2008)\pages
Article 08.3.4, 3pp (electronic)\endref

\ref\key F\by J. Franel\paper On a question of Laisant
\jour L'Interm\'ediaire des Math\'ematiciens\vol 1\yr 1894\pages 45--47\endref

\ref\key G\by H.W. Gould\book {\rm Combinatorial Identities} \publ
Morgantown Printing and Binding Co., 1972\endref

\ref\key JV\by F. Jarvis, H.A. Verrill\paper {\rm Supercongruences for the Catalan-Larcombe-French numbers}
\jour Ramanujan J. 22(2010) 171--186\endref

\ref\key L\by E. Lehmer\paper {\rm On congruences involving Bernoulli numbers and the quotients
of Fermat and Wilson}\jour Ann. of Math. 39(1938) 350--360\endref

\ref\key Sl\by N.J.A. Sloane\paper {\rm Sequence A000172 in OEIS
(On-Line Encyclopedia of Integer Sequences)}
\jour {\tt http://oeis.org/A000172}\endref

\ref\key St92\by V. Strehl\paper {\rm Recurrences and Legendre transform}
\jour S\'em. Lothar. Combin. 29(1992)\pages 1-22\endref

\ref\key St94\by V. Strehl\paper {\rm Binomial identities--combinatorial and algorithmic aspects}
\jour Discrete Math. 136(1994) 309--346\endref

\ref\key S11\by Z.-W. Sun\paper {\rm Super congruences and Euler
numbers} \jour Sci. China Math. 54(2012) 2509-2535\endref

\ref\key S12a\by Z.-W. Sun\paper {\rm Arithmetic theory of harmonic numbers}
\jour Proc. Amer. Math. Soc. 140(2012) 415--428\endref

\ref\key S12b\by Z.-W. Sun\paper {\rm On sums of Ap\'ery polynomials and
related congruences} \jour J. Number Theory 132(2012) 2673--2699\endref

\ref\key S13\by Z.-W. Sun\paper {\rm Connections between $p=x^2+3y^2$ and Franel numbers}
\jour J. Number Theory 133(2013), 2914--2928\endref

\ref\key ST10\by Z.-W. Sun, R. Tauraso\paper {\rm New congruences for central binomial coefficients}
\jour Adv. in Appl. Math. 45(2010) 125--148\endref

\ref\key ST11\by Z.-W. Sun, R. Tauraso\paper {\rm On some new congruences for binomial coefficients}
\jour Int. J. Number Theory 7(2011) 645--662\endref

\ref\key W\by J. Wolstenholme\paper {\rm On certain properties of prime numbers}\jour Quart. J. Appl. Math.
5(1862) 35--39\endref

\ref\key Z\by D. Zagier\paper {\rm Integral solutions of Ap\'ery-like recurrence equations}
\jour in: Groups and Symmetries: from Neolithic Scots to John McKay, CRM Proc. Lecture Notes 47, Amer. Math. Soc.,
Providence, RI, 2009, pp. 349--366\endref


\endRefs
\enddocument